\documentclass[a4paper]{article}

\usepackage{amsthm}
\usepackage{amssymb}
\usepackage{latexsym}
\usepackage{longtable}
\usepackage{epsfig}
\usepackage{amsmath}
\usepackage{hhline}

\newtheorem{thm}{Theorem}[section]
\newtheorem{defn}[thm]{Definition}

\newtheorem{lem}[thm]{Lemma}
\newtheorem{cor}[thm]{Corollary}

\newtheorem{prop}[thm]{Proposition}

\begin{document}

\title{On defect groups for generalized blocks of the symmetric group}
\author{Jean-Baptiste Gramain \\ \\Department of Mathematical Sciences\\University of Aberdeen,
Scotland}
\date{March, 2007}
\maketitle

\begin{abstract}
In a paper of 2003, B. K\"ulshammer, J. B. Olsson and G. R. Robinson
defined $\ell$-blocks for the symmetric groups, where $\ell >1$ is
an arbitrary integer. In this paper, we give a definition for the
defect group of the principal $\ell$-block. We then check that, in
the Abelian case, we have an analogue of one of M. Brou\'e's
conjectures.
\end{abstract}

\section{Introduction}

\subsection{Generalized blocks}

The concept of {\emph{generalized blocks}} was introduced by B.
K\"{u}lshammer, J. B. Olsson and G. R. Robinson in \cite{KOR}. Take
any finite group $G$, and take a union ${\cal C}$ of conjugacy
classes of $G$ containing the identity. We can consider the
restriction to ${\cal C}$ of the ordinary scalar product on
characters of $G$. We denote by Irr$(G)$ the set of complex
irreducible characters of $G$. For $\chi, \, \psi \in
{\mbox{Irr}}(G)$, we let
$$ \langle \chi, \, \psi \rangle_{\cal C}:= \displaystyle \frac{1}{|G|} \sum_{g \in {\cal C}} \chi(g) \psi(g^{-1}).$$
We call $\langle \chi, \, \psi\rangle _{\cal C}$ the
{\emph{contribution}} of $\chi$ and $\psi$ across ${\cal C}$. Then
$\chi$ and $\psi$ are said to be {\emph{directly}} ${\cal
C}${\emph{-linked}} if $\langle \chi, \, \psi\rangle _{\cal C} \neq
0$, and {\emph{orthogonal across}} ${\cal C}$ otherwise. Then direct
${\cal C}$-linking is a reflexive (since $1 \in {\cal C}$) and
symmetric binary relation on Irr$(G)$. Extending it by transitivity,
we obtain an equivalence relation (called ${\cal
C}${\emph{-linking}}) on Irr$(G)$ whose equivalence classes are
called the ${\cal C}$-blocks. Note that, since they are orthogonal
on the whole of $G$, two {\bf{distinct}} irreducible characters are
directly linked across ${\cal C}$ if and only if they are directly
linked across $G \setminus {\cal C}$. Using direct $(G \setminus
{\cal C})$-linking (which needs to be extended by reflexivity and
transitivity), the $(G \setminus {\cal C})$-blocks we obtain are
thus the same as the ${\cal C}$-blocks. Note also that, if we take
${\cal C}$ to be the set of $p$-regular elements of $G$ (i.e. whose
order is not divisible by $p$), for some prime $p$, then the ${\cal
C}$-blocks are just the "ordinary" $p$-blocks (cf e.g.
\cite{Navarro}).

\medskip
In \cite{KOR}, the authors have defined $\ell$-blocks for the
symmetric groups, where $\ell \geq 2$ is any integer. To obtain
this, they take ${\cal C}$ to be the set of $\ell${\emph{-regular}}
elements, i.e. none of whose cycle has length divisible by $\ell$
(in particular, if $\ell$ is a prime $p$, then the $\ell$-blocks are
just the $p$-blocks). The complement of ${\cal C}$ is the set of
$\ell${\emph{-singular}} elements. The $\ell$-blocks thus obtained
satisfy an analogue of the Nakayama Conjecture: two irreducible
characters $\varphi_{\lambda}$ and $\varphi_{\mu}$ of the symmetric
group $S_n$ (where $\lambda$ and $\mu$ are partitions of $n$) belong
to the same $\ell$-block if and only if $\lambda$ and $\mu$ have the
same $\ell$-core.

Following this work, A. Mar\'oti studied generalized blocks in the
alternating groups, and proved that, if $\ell$ is 2 or any odd
integer greater than 1, then the $\ell$-blocks of the alternating
groups also satisfy an analogue of the Nakayama Conjecture (cf
\cite{Attila}).

\subsection{Defect groups}

One key notion in block theory is that of {\emph{defect group}}. The
defect group of a block, or its normalizer, while being smaller
groups than the one we start from, contain a lot of information
about the block and its structure. Particularly striking are
Brou\'e's conjectures, as explained in the next section.

However, we don't have any notion of defect group for generalized
blocks. It is the goal of this paper to give a candidate for the
defect of the principal block in the case of symmetric groups, and
to test its relevance by proving that, in the Abelian case, it
satisfies a weak analogue of one of Brou\'e's conjectures.

\medskip If $G$ is a finite group and $p$ is a prime, we call
{\emph{principal $p$-block}} the $p$-block containing the trivial
character. Its defect group is a Sylow $p$-subgroup of $G$. If we
take $G$ to be the symmetric group $S_n$, and $\ell \geq 2$ to be
any integer, we define the {\emph{principal $\ell$-block}} to be the
$\ell$-block of $G$ containing the trivial character. We would like
to define its defect group to be a Sylow $\ell$-subgroup of $G$, if
any such thing existed. It turns out that, in the symmetric group,
we can construct some kind of $\ell$-analogue of the Sylow
subgroups.

\medskip
Let $p$ be a prime. For this description, we refer to \cite{Olsson}.
We define the $p$-groups $(P_i)_{i \geq 0}$ as follows: let $P_0=1$,
and, for $i>0$, let $P_i= P_{i-1} \wr {\bf Z}_p$ (where ${\bf Z}_p$
is the cyclic group of order $p$). Then, for any $i \geq 0$, $P_i$
is a Sylow $p$-subgroup of $S_{p^i}$. Now take any $n>0$, and let
$a_0+a_1p+ \, \cdots \, +a_kp^k$ be the $p$-adic decomposition of
$n$. Then $P=P_0^{a_0} \times P_1^{a_1} \times \, \cdots \, \times
P_k^{a_k}$ is a Sylow $p$-subgroup of $S_n$. Note that $P$ is
Abelian if and only if $n=a_0+a_1p$ ($0 \leq a_0, a_1 <p$).

\medskip
If we take now any integer $\ell \geq 2$ instead of $p$, we can
still do the same construction. We let $L_0=1$, and $L_i=L_{i-1} \wr
{\bf Z}_{\ell}$ for $i>0$. Then, for any integer $n>0$ with
$\ell$-adic decomposition $a_0+a_1 \ell+ \, \cdots \, +a_k \ell^k$,
we let \\${\cal L}=L_0^{a_0} \times L_1^{a_1} \times \, \cdots \, \times
L_k^{a_k}$. Then ${\cal L}$ is a subgroup of $S_n$, and its order is a
power of $\ell$. We take ${\cal L}$ as {\emph{generalized defect group}}
for the principal $\ell$-block of $S_n$. Note that ${\cal L}$ is Abelian if
and only if $n=a_0+a_1\ell$ ($0 \leq a_0, a_1 <\ell$), and cyclic if
and only if $n=a_0+\ell$ ($0 \leq a_0 <\ell$).

\subsection{Perfect isometries}

The strong link between the principal block and it's defect group is
illustrated by Brou\'e's conjectures (cf \cite{Broue}). One of
Brou\'e's conjectures (which is just the shadow, at the level of
characters, of much deeper equivalences conjectured by Brou\'e)
states that, if $G$ is a finite group with Abelian Sylow
$p$-subgroup $P$, and if $N=N_G(P)$ is the normalizer of $P$ in $G$,
then the principal $p$-blocks of $N$ and $G$ are {\emph{perfectly
isometric}}.

\medskip
This is known to be true when $G$ is a symmetric group (cf
\cite{Raph}). We would like to prove an analogue of this in the case
of generalized blocks for the symmetric group. The tool to do this
is provided in \cite{KOR}. It is the concept of {\emph{generalized
perfect isometry}}, and goes as follows. If $G$ and $H$ are finite
groups, ${\cal C}$ and ${\cal D}$ are closed unions of conjugacy
classes of $G$ and $H$ respectively, and if $b$ (resp. $b'$) is a
union of ${\cal C}$-blocks of $G$ (resp. ${\cal D}$-blocks of $H$),
then we say that there is a {\emph{generalized perfect isometry}}
between $b$ and $b'$ (with respect to ${\cal C}$ and ${\cal D}$) if
there exists a bijection with signs between $b$ and $b'$, which
furthermore preserves contributions; i.e. there exists a bijection
$I \colon b \longmapsto b'$ such that, for each $\chi \in b$, there
is a sign $\varepsilon (\chi )$, and such that
$$\forall \chi, \, \psi \in b, \, \; \langle I(\chi), \, I(\psi) \rangle _{{\cal D}}
=\langle \varepsilon (\chi) \chi, \, \varepsilon (\psi) \psi\rangle
_{{\cal C}}.$$ Note that this is equivalent to $\langle I(\chi), \,
I(\psi) \rangle _{{\cal D'}}=\langle \varepsilon (\chi) \chi, \,
\varepsilon (\psi) \psi\rangle _{{\cal C'}}$, where ${\cal C'}=G
\setminus {\cal C}$ and ${\cal D'}=H \setminus {\cal D}$.

\medskip
When specialized to $p$-blocks, this notion is a bit weaker than the
definition of Brou\'e. If two $p$-blocks $b$ and $b'$ are perfectly
isometric in Brou\'e's sense, then there is a generalized perfect
isometry with respect to $p$-regular elements between $b$ and $b'$.
It is however possible to exhibit generalized perfect isometries in
some cases where there is no perfect isometry in Brou\'e's sense (cf
\cite{JB}).

\medskip
We will prove in Theorem 4.1 that, if $\ell \geq 2$ is any integer,
and if $n=\ell w + r$ for some $0 \leq r, \, w < \ell$, then there
is a generalized perfect isometry with respect to $\ell$-regular (or
$\ell$-singular) elements between the principal $\ell$-blocks of
$S_n$ and $N_{S_n}({\cal L})$, where ${\cal L}$ is defined as in the
previous section (in particular, we only deal with the case where
${\cal L}$ is Abelian). In section 2, we show that the study of the
$\ell$-blocks of $N_{S_n}({\cal L})$ reduces to that of the
$\ell$-blocks of $N_{S_{\ell w}}({\cal L})$. We then compute the
$\ell$-blocks of $N_{S_{\ell}}(L)$, where $L$ is cyclic of order
$\ell$. Section 3 is devoted to the study of the principal
$\ell$-block of $N_{S_{\ell w}}({\cal L})$. This includes general
results about wreath products with symmetric groups. Finally, in
section 4, we state and prove a weak analogue of Brou\'e's Abelian
Defect Conjecture for generalized blocks of the symmetric group.

\section{The Abelian defect case}
Let $\ell \geq 2$ be any integer. From now on, we take $G$ to be the
symmetric group $S_{\ell w + r}$, where $0 \leq r, \, w < \ell$. We
let $L={\bf Z}_{\ell} \hookrightarrow S_{\ell}$, and\\${\cal
L}=<\omega_1> \times \cdots \times <\omega_w> \cong L^w$, where
$\omega_1, \, \ldots, \, \omega_w$ are disjoint $\ell$-cycles of
$G$.

\subsection{Blocks of $N_G({\cal L})$}

We have
$${\cal L} \cong {\bf Z}_{\ell}^w \hookrightarrow S_{\ell w} \hookrightarrow
S_{\ell w} \times S_r \hookrightarrow G=S_{\ell w + r},$$ and
$N_G({\cal L}) \cong N_{S_{\ell w}} ({\cal L}) \times S_r$. Thus
Irr$(N_G({\cal L}))=$Irr$(N_{S_{\ell w }} ({\cal L}))
\otimes$Irr$(S_r)$.

We will define later what we mean by an $\ell${\emph{-singular}}
element of $N_{S_{\ell w }} ({\cal L})$. We define the
$\ell${\emph{-singular}} elements of $N_G({\cal L})\cong N_{S_{\ell
w}} ({\cal L}) \times S_r$ to be the $(\rho, \,
 \sigma)$, where $\rho \in N_{S_{\ell w }} ({\cal L})$ is $\ell$-singular and
$\sigma \in S_r$. For $\psi=\psi_0 \otimes \psi_1$ and $\psi'=\psi'_0
\otimes \psi'_1$ in Irr$(N_G({\cal L}))$, we have
$$ \begin{array}{rl} \langle \psi, \psi'\rangle _{N_G({\cal L}), \, \ell-sing} &= \displaystyle
\frac{1}{|N_G({\cal L})|} \sum_{\scriptstyle g \in N_G({\cal L})
\atop \ell-sing} \psi(g) \overline{\psi'(g)} \\ &=\displaystyle
\frac{1}{|N_G({\cal L})|} \sum_{\scriptstyle \rho \in N_{S_{\ell
w}}({\cal L}) \atop \ell-sing} \sum_{\sigma \in S_r} \psi_0(\rho)
\psi_1(\sigma)
 \overline{\psi'_0(\rho) \psi'_1(\sigma)} \\
&= \delta_{\psi_1 \psi'_1} \langle \psi_0, \psi'_0\rangle
_{N_{S_{\ell w}}({\cal L}), \, \ell-sing}. \end{array}$$

Thus $\psi_0 \otimes \psi_1$ and $\psi'_0 \otimes \psi'_1$ are
directly linked across the set of $\ell$-singular elements of
$N_G({\cal L})$ if and only if $\psi_1=\psi'_1$, and $\psi_0$ and
$\psi'_0$ are directly linked across $\ell$-singular elements of
$N_{S_{\ell w }}({\cal L})$. Hence the $\ell$-blocks of $N_G({\cal
L})$ are the $b \otimes \{ \psi_1 \}$, with $b$ an $\ell$-block of
$N_{S_{\ell w}}({\cal L})$ and $\psi_1 \in$Irr$(S_r)$. In
particular, the principal $\ell$-block of $N_G({\cal L})$ is $b_0
\otimes \{ 1_{S_r} \}$, where $b_0$ is the principal $\ell$-block of
$N_{S_{\ell w}}({\cal L})$.

\subsection{Blocks of $N_{S_{\ell}}(L)$}

We will see later that, in fact, the principal $\ell$-block of
$N_{S_{\ell w}}({\cal L})$ can be obtained from the principal
$\ell$-block of $N_{S_{\ell}}(L)$.

We write $N=N_{S_{\ell}}(L)$. Then $L=<\omega>$ is a regular,
transitive subgroup of $S_{\ell}$, so that, by
\cite{Dixon-Mortimer}, Corollary 4.2 B, we have $N \cong L \rtimes
Aut(L)$, where the action of $Aut(L)$ on $L$ is the natural one.

We have $K:=Aut(L)=\{ \omega \longmapsto \omega^k, \, (k, \, \ell)=1
\}$. In particular, $|K|= \varphi (\ell)$, where $\varphi$ is the
Euler function.

$N=N_{S_{\ell}}(L)$ being a natural subgroup of $S_{\ell}$, we say
that an element of $N_{S_{\ell}}(L)$ is $\ell${\emph{-singular}} if
it is $\ell$-singular as an element of $S_{\ell}$, i.e. just an
$\ell$-cycle. The $\ell$-singular elements of $N$ are precisely the
$\ell$-cycles of $L=<\omega>$, i.e. the $\omega^k$, $(k, \,
\ell)=1$.

\medskip
We let $\sigma=e^{2i\pi/\ell}$ be a primitive $\ell$-th root of
unity. We then have
$$ \mbox{Irr}(L)=\{ \chi_{[m]}, \, m \in {\bf Z} \}= \{ \chi_{[m]},
\, 1 \leq m \leq \ell \},$$ where, for $m \in {\bf Z}$, $\chi_{[m]}$
is uniquely defined by $\chi_{[m]}(\omega)=\sigma^m$.

We want to use Clifford's Theory to describe the irreducible
characters of $N$. $N$ acts on Irr$(L)$ via $\chi \mapsto \chi^g$,
where, for $g \in N$, we have $\chi^g(h)=\chi(g^{-1}hg)$ for all $h
\in L$. For $m \in {\bf Z}$, we write $I_{[m]}= \{ g \in N \, | \,
\chi^g=\chi \}$ the inertia subgroup of $\chi_{[m]}$, and
$\chi_{[m]}^N=\{ \chi_{[m]}^g, \, g \in N \}$ the orbit of
$\chi_{[m]}$ under the action of $N$. Then Clifford's Theorem has
the following consequence:

\begin{thm}(\cite{Isaacs}, Theorem (6.11))
Take any $ 1 \leq m \leq \ell$. Using the notations above, let
$${\cal A}_{[m]}= \{ \varphi \in \mbox{Irr}(I_{[m]}) \, | \,
\langle \varphi  \mathord{\downarrow}_L, \chi_{[m]}\rangle  \neq 0
\}$$ and
$${\cal B}_{[m]}= \{ \psi \in \mbox{Irr}(N) \, | \,
\langle \psi \mathord{\downarrow} _L, \chi_{[m]}\rangle  \neq 0
\}.$$ Then $\varphi \longmapsto \varphi \mathord{\uparrow}^N$ is a
bijection from ${\cal A}_{[m]}$ to ${\cal B}_{[m]}$ such that
$$\langle (\varphi \mathord{\uparrow}^N )\mathord{\downarrow}_L, \chi_{[m]}\rangle =\langle \varphi \mathord{\downarrow}_L,
\chi_{[m]}\rangle .$$
\end{thm}

\noindent
Now, for any $1 \leq m \leq \ell$, we have
$$ L \lhd I_{[m]} \leq L \rtimes Aut(L) = L \rtimes K,$$
so that $I_{[m]}=L \rtimes K_{[m]}$ for some $K_{[m]} \leq K$. This
implies easily that any irreducible character of $L$ extends to its
inertia subgroup. Hence, by a result of Gallagher (cf \cite{Isaacs},
Corollary (6.17)), we have, for any $1 \leq m \leq \ell$,
$${\cal A}_{[m]} = \{ \tilde{\chi}_{[m]} \vartheta, \, \vartheta \in
\mbox{Irr}(I_{[m]}/L) \} = \{ \tilde{\chi}_{[m]} \vartheta, \,
\vartheta \in \mbox{Irr}(K_{[m]}) \},$$ where $\tilde{\chi}_{[m]}$
is the (unique) extension of $\chi_{[m]}$ to $I_{[m]}$.

\medskip

By Clifford's Theory, for any $\psi \in$Irr$(N)$, the irreducible
components of $\psi \mathord{\downarrow}_L$ are the elements of
exactly one orbit of Irr$(L)$ under the action of $N$. In order to
obtain a description of Irr$(N)$, it therefore suffices to find a
partition of Irr$(L)$ into $N$-orbits.

\medskip
First note that $g \in S_{\ell}$ belongs to $N=
N_{S_{\ell}}(<\omega>)$ if and only if $g^{-1}\omega g=\omega^k$ for
some $(k, \, \ell)=1$. For such a $g$, and for any $1 \leq m \leq
\ell$, we have
$$
\chi_{[m]}^g(\omega)=\chi_{[m]}(\omega^g)=\chi_{[m]}(\omega^k)=\sigma^{km}=\chi_{[km]}(\omega).$$
Hence $\chi_{[m]}^N=\{\chi_{[km]}, \,(k, \, \ell)=1 \}$. In
particular, we have $$\chi_{[m]}^N \subset \{\chi_{[n]}, \,(n, \,
\ell)=(m, \, \ell) \}.$$ On the other hand, we have
$$ \mbox{Irr}(L)= \displaystyle \coprod_{d|\ell} \{\chi_{[n]}, \,(n, \, \ell)=d
\}.$$We partition Irr$(L)$ as follows:

\begin{prop}
We have
$$ \mbox{Irr}(L)= \displaystyle \coprod_{d|\ell} \chi_{[d]}^N,$$
and, for all $d| \ell$, we have $| \chi_{[d]}^N | =
\varphi(\frac{\ell}{d})$.
\end{prop}

\begin{proof}

We start by noting that, for $d | \ell$ and $1 \leq n \leq \ell$, we
have $(n, \, \ell)=d$ if and only if $n=dk$ for some $1 \leq k \leq
\frac{\ell}{d}$ and $(k, \, \frac{\ell}{d})=1$. Thus
$$\{\chi_{[n]}, \,(n, \, \ell)=d
\}=\{\chi_{[kd]}, \, 1 \leq k \leq \frac{\ell}{d} , \,(k, \,
\frac{\ell}{d})=1 \},$$ and$$ \mbox{Irr}(L)= \displaystyle
\coprod_{d|\ell} \{\chi_{[n]}, \,(n, \, \ell)=d \}=\displaystyle
\coprod_{d|\ell}\{\chi_{[kd]}, \, 1 \leq k \leq \frac{\ell}{d} ,
\,(k, \, \frac{\ell}{d})=1 \}.$$ Write $\ell = \prod_{i \in {\cal
I}} p_i^{a_i}$, where $a_i
>0$ for $i \in {\cal I}$, and $p_i \neq p_j$ for $i \neq j$. Now
take $d| \ell$, and write $d = \prod_{i \in {\cal I}_d} p_i^{b_i}$,
where $0< b_i \leq a_i$ for $i \in {\cal I}_d \subset {\cal I}$.
Then
$$\displaystyle \frac{\ell}{d} =\prod_{i \in {\cal I}_d} p_i^{a_i-b_i} \prod_{i \in {\cal I}
\setminus {\cal I}_d} p_i^{a_i}.$$

Suppose first that $b_i < a_i$ for all $i \in {\cal I}_d$. Then, for
any $k \in {\bf N}$, $(k, \, \frac{\ell}{d})=1$ implies $(k, \,
d)=1$. Thus, in this case, for any $1 \leq k \leq \frac{\ell}{d}$
such that $(k, \, \frac{\ell}{d})=1$, we have $(k, \, \ell)=1$, so
that $\chi_{[kd]} \in \chi_{[d]}^N$. Hence
$$\chi_{[d]}^N=\{\chi_{[kd]}, \,(k, \, \ell)=1 \}=
\{\chi_{[kd]}, \, 1 \leq k \leq \frac{\ell}{d} , \,(k, \,
\frac{\ell}{d})=1 \}.$$ In particular, $| \chi_{[d]}^N | =
\varphi(\frac{\ell}{d})$.

\medskip
Suppose now that there exists $i \in {\cal I}_d$ such that
$b_i=a_i$. Take $1 \leq k \leq \frac{\ell}{d}$ such that $(k, \,
\frac{\ell}{d})=1$. Write $k=k_dk_{\ell}$, where $(k_{\ell}, \,
\ell)=1$, and $k_d=\prod_{i \in {\cal I}_k} p_i^{c_i}$, with ${\cal
I}_k \subset {\cal I}_d$ and $c_i >0$ for $i \in {\cal I}_k$. Let
$m=\prod_{i \in {\cal I}_d \setminus {\cal I}_k} p_i$ and $n=m
\frac{\ell}{d} +k_d$. Now take any $i \in {\cal I}$.
\begin{itemize}
\item{}
If $i \in {\cal I}_k$, then $p_i | k_d$, and $(k_d, \,
\frac{\ell}{d})=1$, so that $p_i  \not | \, \frac{\ell}{d}$, and
$p_i \not | \, m$. Thus $p_i \not | \, n$.
\item{}
If $i \in {\cal I} \setminus {\cal I}_k$, then $p_i | m
\frac{\ell}{d}$, and $p_i \not | \, k_d$. Thus $p_i \not | \, n$.
\end{itemize}
Hence $(n, \, \ell)=1$ and, since $(k_{\ell}, \, \ell)=1$, we have
$(nk_{\ell}, \, \ell)=1$. Now
$$ \begin{array}{rl} k-nk_{\ell} & =k_dk_{\ell}-nk_{\ell} \\ &
=k_{\ell}(k_d-(m \frac{\ell}{d} +k_d)) \\ & =-k_{\ell}m
\frac{\ell}{d} \equiv 0 \; (\frac{\ell}{d}). \end{array}$$Thus $k
\equiv n k_{\ell} \; (\frac{\ell}{d})$ and $kd \equiv n k_{\ell}d \;
(\ell)$, whence $\chi_{[kd]}=\chi_{[nk_{\ell}d]}$ and $(nk_{\ell},
\, \ell)=1$. Hence $\chi_{[kd]} \in \chi_{[d]}^N$. Thus, in this
case too,
$$\chi_{[d]}^N=\{\chi_{[kd]}, \,(k, \, \ell)=1 \}=
\{\chi_{[kd]}, \, 1 \leq k \leq \frac{\ell}{d} , \,(k, \,
\frac{\ell}{d})=1 \}.$$ In particular, $| \chi_{[d]}^N | =
\varphi(\frac{\ell}{d})$.

Finally, we get
$$ \mbox{Irr}(L)=\displaystyle \coprod_{d|\ell}\{\chi_{[kd]}, \,
1 \leq k \leq \frac{\ell}{d} , \,(k, \, \frac{\ell}{d})=1 \}=
\displaystyle \coprod_{d|\ell} \chi_{[d]}^N.$$
\end{proof}

This gives us the whole of Irr$(N_{S_{\ell}}(L))$:
$$ \begin{array}{rl} \mbox{Irr}(N_{S_{\ell}}(L)) & = \displaystyle \coprod_{d | \ell}
\left\{ (\tilde{\chi}_{[d]} \vartheta)
\mathord{\uparrow}^{N_{S_{\ell}}(L)}, \, \vartheta \in
\mbox{Irr}(I_{[d]}/L) \right\} \\ & = \{ \psi_{d,\vartheta}, \, d|
\ell, \, \vartheta \in \mbox{Irr}(I_{[d]}/L) \}. \end{array}$$Note
that, for $d | \ell$, we have $[N_{S_{\ell}}(L) \colon I_{[d]}]=|
\chi_{[d]}^N | = \varphi(\frac{\ell}{d})$, so that
$$|I_{[d]}/L| = \displaystyle \frac{[N_{S_{\ell}}(L) \colon L]}{[N_{S_{\ell}}(L) \colon
I_{[d]}]}=\frac{\varphi(\ell)}{\varphi(\frac{\ell}{d})}.$$Furthermore,
for any $d | \ell$ and $\vartheta \in \mbox{Irr}(I_{[d]}/L)$, we
have
$$\psi_{d,\vartheta}\mathord{\downarrow}_{L} = \displaystyle \sum_{ \scriptstyle 1 \leq n \leq \ell
 \atop
\scriptstyle (n , \, \ell)=d }  \chi_{[n]} = \sum_{ \scriptstyle 1
\leq k \leq \frac{\ell}{d} \atop \scriptstyle (k, \,
\frac{\ell}{d})=1 } \chi_{[kd]}.$$

\medskip
We can now construct the $\ell$-blocks of $N$. For any $d,d' |
\ell$, $\vartheta \in \mbox{Irr}(I_{[d]}/L)$ and $\vartheta' \in
\mbox{Irr}(I_{[d']}/L)$, we have
$$\langle \psi_{d,\vartheta},\psi_{d',\vartheta'}\rangle _{\ell-sing}= \displaystyle
\frac{1}{\ell \varphi(\ell)} \sum_{ \scriptstyle 1 \leq k \leq \ell
\atop \scriptstyle (k, \, \ell)=1 }  \psi_{d,\vartheta} (\omega^k)
 \overline{\psi_{d',\vartheta'}(\omega^k)} = \frac{1}{\ell} \psi_{d,\vartheta} (\omega)
 \overline{\psi_{d',\vartheta'}(\omega)}$$
(all $\ell$-cycles of $L$ being conjugate in $N$). Now, by the
above, we have
$$\psi_{d,\vartheta} (\omega) = \displaystyle \sum_{ \scriptstyle 1 \leq n \leq \ell
 \atop
\scriptstyle (n , \, \ell)=d }  \chi_{[n]}(\omega)=\sum_{
\scriptstyle 1 \leq n \leq \ell
 \atop
\scriptstyle (n , \, \ell)=d } \sigma^n$$(and similarly for
$\psi_{d',\vartheta'}(\omega)$).

\begin{lem}
For any $d | \ell$, $\{ \sigma^n, \, 1 \leq n \leq \ell , \, (n, \,
\ell)=d \}$ is the set of primitive $\frac{\ell}{d}$-th roots of
unity.
\end{lem}
\begin{proof}
Take any $1 \leq n \leq \ell$ with $(n, \, \ell)=d$. Write $n=kd$,
for some $1 \leq k \leq \frac{\ell}{d}$, $(k, \, \frac{\ell}{d})=1$.
Then $\sigma^{n \frac{\ell}{d}}=\sigma^{k\ell}=1$, so that
$\sigma^n$ is an $\frac{\ell}{d}$-th root of 1.

Suppose $\sigma^{nm}=1$. Then $\sigma^{kdm}=1$. Thus $\ell | kdm$,
so that $\frac{\ell}{d} | km$, and $\frac{\ell}{d} | m$ (since $(k,
\, \frac{\ell}{d})=1$). Hence $\sigma^n$ is a primitive
$\frac{\ell}{d}$-th root of 1.

To conclude, we just note that $1 \leq n \neq n' \leq \ell$ implies
$\sigma^n \neq \sigma^{n'}$, so that
$$|\{ \sigma^n, \, 1 \leq n \leq \ell , \, (n, \,
\ell)=d \}|=|\{ 1 \leq n \leq \ell , \, (n, \, \ell)=d \}|=
\varphi(\frac{\ell}{d}),$$the number of primitive
$\frac{\ell}{d}$-th roots of 1.
\end{proof}

Thus, for any $d | \ell$ and $\vartheta \in \mbox{Irr}(I_{[d]}/L)$,
$\psi_{d,\vartheta} (\omega)$ is the sum of the primitive
$\frac{\ell}{d}$-th roots of unity. This is an example of
{\emph{Ramanujan's sum}} (cf \cite{HW} 5.6 (2)). Ramanujan's sum
$c_q(m)$, defined for integers $q>0$ and $m \geq 0$, is the sum of
the $m$-th powers of the primitive $q$-th roots of unity. We thus
have $\psi_{d,\vartheta} (\omega)=c_{\frac{\ell}{d}}(1)$. Now, by
\cite{HW} Theorem 271, we have, for any $q>0$ and $m \geq 0$,
$$c_q(m)= \displaystyle \sum_{ \scriptstyle n | m, \, n| q} \mu
\left( \frac{q}{n} \right) n,$$where $\mu$ is the M\"obius function
(cf \cite{HW} 16.3). In particular,
$c_{\frac{\ell}{d}}(1)=\mu(\frac{\ell}{d})$.

Finally, we get that, for any $d,d' | \ell$, $\vartheta \in
\mbox{Irr}(I_{[d]}/L)$ and $\vartheta' \in \mbox{Irr}(I_{[d']}/L)$,
$$\langle \psi_{d,\vartheta},\psi_{d',\vartheta'}\rangle _{\ell-sing}= \displaystyle
\frac{1}{\ell} \; \mu \left( \frac{\ell}{d} \right) \mu \left(
\frac{\ell}{d'} \right).$$This allows us to compute the
$\ell$-blocks of $N_{S_{\ell}}(L)$. Take any $d|\ell$.
\begin{itemize}
\item{}
If $\mu(\frac{\ell}{d})=0$, then $\langle
\psi_{d,\vartheta},\psi_{d',\vartheta'}\rangle _{\ell-sing}= 0$ for
any $d' | \ell$, $\vartheta \in \mbox{Irr}(I_{[d]}/L)$ and
$\vartheta' \in \mbox{Irr}(I_{[d']}/L)$. Thus
$\{\psi_{d,\vartheta}\}$ is an
$\ell$-block of $N_{S_{\ell}}(L)$ for any\\
$\vartheta \in
\mbox{Irr}(I_{[d]}/L)$.
\item{}
If $\mu(\frac{\ell}{d}) \neq 0$, then $\langle
\psi_{d,\vartheta},\psi_{d',\vartheta'}\rangle _{\ell-sing} \neq 0$
for any $d' | \ell$ such that $\mu(\frac{\ell}{d'}) \neq 0$, and any
$\vartheta \in \mbox{Irr}(I_{[d]}/L)$ and $\vartheta' \in
\mbox{Irr}(I_{[d']}/L)$. Thus
$$\{\psi_{d,\vartheta} \, ; \; d | \ell, \, \mu(\frac{\ell}{d}) \neq 0 , \, \vartheta
\in \mbox{Irr}(I_{[d]}/L) \}$$is an $\ell$-block of
$N_{S_{\ell}}(L)$, namely the principal block $\beta_0$ since it
contains the trivial character $1_{N_{S_{\ell}}(L)}=\psi_{\ell,1}$.
\end{itemize}

Note in particular that $N_{S_{\ell}}(L)$ has a single $\ell$-block
if and only if $\ell$ is squarefree.

\medskip
We now want to compute the order of $\beta_0$. We have
$$ \begin{array}{rl} |\beta_0| & = \displaystyle \sum_{ \scriptstyle d|
\ell, \, \mu(\frac{\ell}{d}) \neq 0} | \mbox{Irr}(I_{[d]}/L) | \\
 & = \displaystyle \sum_{ \scriptstyle d|
\ell, \, \mu(\frac{\ell}{d}) \neq 0} | I_{[d]}/L |  \; \; \; \;
(\mbox{since} \; I_{[d]}/L \leq Aut(L) \; \mbox{is Abelian})\\  & =
\displaystyle \sum_{ \scriptstyle d| \ell, \, \mu(\frac{\ell}{d})
\neq 0} \frac{\varphi(\ell)}{\varphi(\frac{\ell}{d})} \; \; \; \;
\mbox{(as we remarked before).} \end{array}$$

Now, for $d | \ell$, we have $\mu (\frac{\ell}{d}) \neq 0$ if and
only if $\frac{\ell}{d}$ is squarefree. Write $\ell= \prod_{i \in
{\cal I}} p_i^{a_i}$, where the $p_i$'s are distinct primes and $a_i
>0$ for $i \in {\cal I}$. Let $\ell_0= \prod_{i \in {\cal I}} p_i$.
For $d | \ell$, we have $\mu (\frac{\ell}{d}) \neq 0$ if and only if
$\frac{\ell}{d}= \prod_{i \in {\cal J}} p_i$ for some ${\cal J}
\subset {\cal I}$ (in particular, $\mu (\frac{\ell}{d}) \neq 0$ if
and only if $\frac{\ell}{d} | \ell_0$). Now, by \cite{HW} Theorem
62, we have $\varphi(\ell) = \ell \prod_{i \in {\cal I}} (1 -
\frac{1}{p_i})$, and, if $\frac{\ell}{d}= \prod_{i \in {\cal J}}
p_i$, then $\varphi(\frac{\ell}{d}) = \frac{\ell}{d} \prod_{i \in
{\cal J}} (1 - \frac{1}{p_i})$. Thus
$\frac{\varphi(\ell)}{\varphi(\frac{\ell}{d})} = d \prod_{i \in
{\cal I} \setminus {\cal J}} (1 - \frac{1}{p_i})$.

Now $\frac{\ell_0 d}{\ell}=\prod_{i \in {\cal I} \setminus {\cal J}}
p_i$. Thus $\prod_{i \in {\cal I} \setminus {\cal J}} (1 -
\frac{1}{p_i})= \varphi( \frac{\ell_0 d}{\ell}) \frac{\ell}{\ell_0
d}$, \\
and $\frac{\varphi(\ell)}{\varphi(\frac{\ell}{d})} =
\frac{\ell}{\ell_0} \varphi( \frac{\ell_0 d}{\ell})$. We get
$$ |\beta_0|=\displaystyle \sum_{ \scriptstyle d| \ell, \, \mu(\frac{\ell}{d}) \neq 0}
\frac{\varphi(\ell)}{\varphi(\frac{\ell}{d})}=\sum_{\scriptstyle
\delta | \ell_0} \frac{\ell}{\ell_0} \varphi \left(
\frac{\ell_0}{\delta} \right) = \frac{\ell}{\ell_0}
\sum_{\scriptstyle \delta | \ell_0} \varphi \left(
\frac{\ell_0}{\delta} \right) = \ell.$$

\section{Principal block of $N_{S_{\ell w}}({\cal L})$}

We now need to compute the principal $\ell$-block of ${\cal
N}=N_{S_{\ell w}}({\cal L})$. We have ${\cal N}=N_{S_{\ell w}}( L^w)
\cong N_{S_{\ell}}(L) \wr S_w=N \wr S_w$. In order to do this study,
we need some results about the conjugacy classes and irreducible
characters of wreath products.

\subsection{Wreath products}

Let $H$ be a finite group (we will later use the following general
results for $H=L$ and $H=N$). The following facts can be found in
\cite{Pfeiffer}.

Take $g_1, \, \ldots , \, g_r$ representatives for the conjugacy
classes of $H$. The elements of the wreath product $H \wr S_w$ are
of the form $(h;\gamma)=(h_1, \, \ldots, \, h_w ; \gamma)$, with
$h_1, \, \ldots, \, h_w  \in H$ and $ \gamma \in S_w$. For any such
element, and for any $k$-cycle \\$\kappa=(j, j\kappa, \ldots, j
\kappa^{k-1})$ in $\gamma$, we define the {\emph{cycle product}} of $(h;\gamma)$
and $\kappa$ by
$$g((h;\gamma),\kappa)= h_j h_{j \kappa^{-1}} h_{j \kappa^{-2}} \ldots h_{j
\kappa^{-(k-1)}}.$$In particular, $g((h;\gamma),\kappa) \in H$. If
$\gamma$ has cycle structure $\pi$ say, then we form $r$ partitions
$(\pi_1, \, \ldots, \, \pi_r)$ from $\pi$ as follows: any cycle
$\kappa$ in $\pi$ gives a cycle of the same length in $\pi_i$ if the
cycle product $g((h;\gamma),\kappa)$ is conjugate to $g_i$. The
resulting $r${\emph{-tuple of partitions}} of $w$ describes the
{\emph{cycle structure}} of $(h;\gamma)$. We then have

\begin{lem}(\cite{Pfeiffer} Proposition 4.1)
Two elements of $H \wr S_w$ are conjugate if and only if they have
the same cycle structure, and the conjugacy classes of $H \wr S_w$
are parametrized by the $r$-tuples of partitions of $w$: $$\{ (\pi_1,
\, \ldots, \, \pi_r) \, ; \displaystyle \sum_{i=1}^r | \pi_i | =w \}=\{ (\pi_1, \,
\ldots, \, \pi_r) \Vdash w \} .$$
\end{lem}

\medskip
Note that, if $H$ is a natural subgroup of $S_{\ell}$, and if we imbed $H \wr S_w$ in $S_{\ell w}$, then, for any $k$-cycle $\kappa$ of $\gamma$, each $m$-cycle of the cycle product $g((h;\gamma),\kappa)$ (seen as an element of $S_{\ell}$) corresponds to $gcd(m,k)$ cycles of length $lcm(m,k)$ in $(h;\gamma)$, when seen as an element of $S_{\ell w}$.

\medskip
We have the following result about centralizers:

\begin{lem}(\cite{Pfeiffer}, Lemma 4.2)
If $g \in H \wr S_w$ has cycle type $(\pi_1, \, \ldots, \, \pi_r)
\Vdash w$, and if $a_{ik}$ denotes the number of $k$-cycles of
$\pi_i$ ($1 \leq i \leq r$, $1 \leq k \leq w$), then
$$|C_{H \wr S_w}(g)| = \displaystyle \prod_{i,k} a_{ik} ! (k
|C_H(g_i)|)^{a_{ik}}.$$
\end{lem}

\bigskip

The irreducible (complex) characters of $H \wr S_w$ are also
canonically labeled by the $r$-tuples of partitions of $w$. The
construction goes as follows:

We write Irr$(H)=\{ \psi_1, \, \ldots , \, \psi_r \}$ and, for any
$m>0$, \\Irr$(S_m)=\{ \varphi_{\lambda}, \, \lambda \vdash m \}$.
For any $\alpha=(\alpha^1, \, \ldots , \, \alpha^r) \Vdash w$, the
irreducible character $\prod_{i=1}^r \psi_i^{|\alpha^i|}$ of the
base group $N^w$ can be extended in a natural way to its inertia
subgroup $H \wr S_{|\alpha^1|} \times \cdots \times H \wr
S_{|\alpha^r|}$, giving the irreducible character $\prod_{i=1}^r
\widehat{\psi_i^{|\alpha^i|}}$. Any irreducible character of $H \wr
S_{|\alpha^1|} \times \cdots \times H \wr S_{|\alpha^r|}$ which
extends $\prod_{i=1}^r \psi_i^{|\alpha^i|}$ is a tensor product of
$\prod_{i=1}^r \widehat{\psi_i^{|\alpha^i|}}$ with an irreducible
character of the inertia factor $S_{|\alpha^1|} \times \cdots \times
S_{|\alpha^r|}$. Inducing this tensor product to $H \wr S_w$ gives
an irreducible character, and any irreducible character of $H \wr
S_w$ can be obtained in this way. We have the following

\begin{lem}(\cite{Pfeiffer}, Proposition 4.3)
The irreducible characters of $H \wr S_w$ are labeled by the
$r$-tuples of partitions of $w$. We have Irr$(H \wr S_w)=\{
\chi^{\alpha} , \, \alpha \Vdash w \}$, where, for
$\alpha=(\alpha^1, \, \ldots , \, \alpha^r) \Vdash w$,
$$\chi^{\alpha}= \displaystyle  \left( \prod_{i=1}^r
\widehat{\psi_i^{|\alpha^i|}} \otimes \varphi_{\alpha^i} \right)^{H
\wr S_w}.$$
\end{lem}

There is a generalization of the Murnaghan-Nakayama Rule in $H \wr
S_w$ which enlightens the link between the parametrizations for
conjugacy classes and irreducible characters.

\begin{thm}(\cite{Pfeiffer}, Theorem 4.4)
Take $\alpha=(\alpha^1, \, \ldots , \, \alpha^r) \Vdash w$. Take \\$g
\in H \wr S_w$ of cycle type $(a_{ij}(g))_{1 \leq i \leq r, 1 \leq j
\leq w}$ such that $a_{tk}(g) >0$ for some $1 \leq t \leq r$ and $1
\leq k \leq w$. Let $\rho \in H \wr S_{w-k}$ be of cycle type
$(a_{ij}(\rho))_{1 \leq i \leq r, 1 \leq j \leq w}$, where
$$(a_{ij}(\rho)) = \left\{ \begin{array}{ll} a_{tk}(g)-1 &
{\mbox{if}} \; i=t \; \mbox{and} \; j=k \\ a_{ij}(g) &
\mbox{otherwise} \end{array} \right. . $$ If we let
$\chi^{\emptyset}=1$, we then have
$$\chi^{\alpha}(g)= \displaystyle \sum_{s=1}^r \psi_s(g_t)
\sum_{h_{ij}^{\alpha^s}=k} (-1)^{L_{ij}^{\alpha^s}}
\chi^{\alpha-R_{ij}^{\alpha^s}}(\rho),$$ where the second sum is
taken over the set of hooks $h_{ij}^{\alpha^s}$ of length $k$ in
$\alpha^s$, $L_{ij}^{\alpha^s}$ is the leg-length of
$h_{ij}^{\alpha^s}$, and $\alpha-R_{ij}^{\alpha^s}$ is the $r$-tuple
of partitions of $w-k$ obtained from $\alpha$ by removing
$h_{ij}^{\alpha^s}$ from $\alpha^s$.
\end{thm}

We will also write the previous result as
$$\chi^{\alpha}(g)= \displaystyle \sum_{s=1}^r \psi_s(g_t)
\sum_{\beta^s \in {\cal L}_{\alpha^s}^k} (-1)^{L_{\alpha^s \beta^s}}
\chi^{\alpha_{s,\beta^s}}(\rho),$$ where ${\cal L}_{\alpha^s}^k$ is
the set of partitions $\beta^s$ which can be obtained from
$\alpha^s$ by removing a $k$-hook (including multiplicities),
$L_{\alpha^s \beta^s}$ is the leg-length of the hook to remove from
$\alpha^s$ to get to $\beta^s$, and $\alpha_{s,\beta^s}=(\alpha^1,
\, \ldots , \, \alpha^{s-1}, \, \beta^s , \, \alpha^{s+1}, \, \ldots
, \, \alpha^r) \Vdash w-k$.

\subsection{The wreath product $L \wr S_w$; a generalized perfect
isometry}

By the results of the previous section, the conjugacy classes of $L
\wr S_w$ are parametrized by the $\ell$-tuples $(\pi_1, \, \ldots ,
\, \pi_{\ell})$ of partitions of $w$. We take representatives $(g_1,
\, \ldots, \, g_{\ell}=1)$ for the conjugacy classes of $L$. For later coherence, we take $g_i=\omega^i$ ($1 \leq i \leq \ell$), where $L=< \omega>$. We
partition $L \wr S_w$ into unions of conjugacy classes according to
$\pi_{\ell} \in {\cal P}_{\leq w}$ (i.e. $\pi_{\ell}$ partition of
at most $w$). We have
$$L \wr S_w = \displaystyle \coprod_{\pi_{\ell} \in {\cal P}_{\leq
w}} {\cal D}_{\pi_{\ell}},$$ where ${\cal D}_{\pi_{\ell}}$ is the
set of elements whose cycle type has $\pi_{\ell}$ as $\ell$-th part.
Then ${\cal D}_{\emptyset}$ is the set of {\emph{regular}} elements
of $L \wr S_w$. The elements of $L \wr S_w \setminus {\cal
D}_{\emptyset}$ are called {\emph{singular}}.

\bigskip
We now turn to the principal $\ell$-block of $S_{\ell w +r}$. The
trivial character of $S_{\ell w+r}$ is labeled by the partition
$(\ell+r)$ of $\ell+r$, whose $\ell$-core is $(r) \vdash r$.
By the Nakayama Conjecture for generalized blocks (\cite{KOR},
Theorem 5.13), the characters of the principal $\ell$-block $B_0$ of
$S_{\ell+r}$ are labeled by the partitions of $\ell+r$ with the same
$\ell$-core $(r)$. We write ${\cal P}_r$ the set of partitions of $\ell+r$ with $\ell$-core $(r)$, and write $B_0=\{ \varphi_{\lambda}, \, \lambda
 \in {\cal P}_r \}$.
The characters of $B_0$ are parametrized by the $\ell$-quotients of
the partitions labeling them. For $\lambda \in {\cal P}_r$, we write
$\alpha_{\lambda}= ( \alpha_{\lambda}^1, \, \ldots , \,
\alpha_{\lambda}^{\ell} ) \Vdash w$ the $\ell$-quotient of
$\lambda$. These $\ell$-quotients also parametrize the irreducible
(complex) characters of $L \wr S_w$. A key result of \cite{KOR} is
the following:

\begin{thm}(\cite{KOR}, Theorem 5.10 and Proposition 5.11)
With the above notations, the map $\lambda \longmapsto
\alpha_{\lambda}$ induces a generalized perfect isometry between
$B_0$ and Irr$(L \wr S_w)$ with respect to $\ell$-regular and
regular elements respectively. That is, there exist signs $\{
\varepsilon(\lambda) , \, \lambda \in {\cal P}_r \}$ such that, for any $\lambda, \mu
\in {\cal P}_r$,
$$\langle \varphi_{\lambda}, \, \varphi_{\mu} \rangle _{\scriptstyle
\ell-reg \atop S_{\ell w +r}} =\langle \varepsilon (\lambda)
\chi^{\alpha_{\lambda}}, \, \varepsilon (\mu)
\chi^{\alpha_{\mu}}\rangle _{\scriptstyle reg \atop L \wr S_w}.$$
\end{thm}

\bigskip
Given two $\ell$-tuples $\alpha$ and $\beta$ of partitions of $w$,
we will give another expression for $\langle \chi^{\alpha}, \,
\chi^{\beta}\rangle _{\scriptstyle reg \atop L \wr S_w}$. We have
$$\langle \chi^{\alpha}, \,
\chi^{\beta}\rangle _{\scriptstyle reg \atop L \wr
S_w}=\delta_{\alpha \beta} - \langle \chi^{\alpha}, \,
\chi^{\beta}\rangle _{\scriptstyle sing \atop L \wr S_w},$$ and,
with the above notations, we have
$$\langle \chi^{\alpha}, \,
\chi^{\beta}\rangle _{\scriptstyle sing \atop L \wr S_w}=
\displaystyle \sum_{ \emptyset \neq \pi_{\ell} \in {\cal P}_{\leq
w}} \langle \chi^{\alpha}, \, \chi^{\beta}\rangle _{\scriptstyle
{\cal D}_{\pi_{\ell}}}.$$ We have, for any $\emptyset \neq
\pi_{\ell} \in {\cal P}_{\leq w}$,
$$\langle \chi^{\alpha}, \, \chi^{\beta}\rangle _{\scriptstyle
{\cal D}_{\pi_{\ell}}}= \displaystyle \frac{1}{|L \wr S_w|}
\sum_{\scriptstyle (\pi_1, \ldots, \pi_{\ell-1},\emptyset) \Vdash
w-|\pi_{\ell}|} |{\cal D}_{(\pi_1, \ldots, \pi_{\ell})}|
\chi^{\alpha}({\cal D}_{(\pi_1, \ldots, \pi_{\ell})} )
\overline{\chi^{\beta}( {\cal D}_{(\pi_1, \ldots, \pi_{\ell})})},$$
where ${\cal D}_{(\pi_1, \ldots, \pi_{\ell})}$ is the conjugacy
class of cycle type $(\pi_1, \ldots, \pi_{\ell})$ of $L \wr S_w$.

Now, for any $(\pi_1, \ldots, \pi_{\ell-1},\emptyset) \Vdash
w-|\pi_{\ell}|$, Lemma 3.2 gives, writing $b_{jk}$ for the number of $k$-cycles in $\pi_j$
($1\leq j \leq \ell$, $1 \leq k \leq w$),
$$\begin{array}{rl} \displaystyle \frac{|L \wr
S_w|}{|{\cal D}_{(\pi_1, \ldots, \pi_{\ell})}|} & = \displaystyle \prod_{\scriptstyle 1\leq j \leq \ell \atop 1 \leq k \leq w} b_{jk}! (k |C_L(g_j)|)^{b_{jk}} \\ & =  \displaystyle \prod_{\scriptstyle 1 \leq k \leq w} b_{\ell k}! (k |C_L(g_{\ell})|)^{b_{\ell k}} \displaystyle \prod_{\scriptstyle 1\leq j \leq \ell-1 \atop 1 \leq k \leq w} b_{jk}! (k |C_L(g_j)|)^{b_{jk}} \\ & = \displaystyle \prod_{\scriptstyle 1 \leq k \leq w} b_{\ell k}! (k \ell)^{b_{\ell k}} \displaystyle \frac{|L \wr
S_{w-|\pi_{\ell}|}|}{|{\cal D}_{(\pi_1, \ldots, \pi_{\ell-1}, \emptyset)}^{w- |\pi_{\ell}|}|}, \end{array}
$$
where ${\cal D}_{(\pi_1, \ldots, \pi_{\ell-1}, \emptyset)}^{w- |\pi_{\ell}|}$ denotes the set of elements of cycle type $(\pi_1, \ldots, \pi_{\ell-1}, \emptyset)$ in $L \wr
S_{w-|\pi_{\ell}|}$. Writing $d_{\pi_{\ell}}=\prod_{\scriptstyle 1 \leq k \leq w} b_{\ell k}! (k \ell)^{b_{\ell k}}$, we get

$$\langle \chi^{\alpha}, \, \chi^{\beta}\rangle _{\scriptstyle
{\cal D}_{\pi_{\ell}}}= \displaystyle \frac{1}{d_{\pi_{\ell}}}
\sum_{\scriptstyle (\pi_1, \ldots, \pi_{\ell-1},\emptyset) \Vdash
w-|\pi_{\ell}|}  \frac{|{\cal D}_{(\pi_1, \ldots, \pi_{\ell-1}, \emptyset)}^{w- |\pi_{\ell}|}|}{|L \wr
S_{w-|\pi_{\ell}|}|}
\chi^{\alpha}({\cal D}_{(\pi_1, \ldots, \pi_{\ell})} )
\overline{\chi^{\beta}( {\cal D}_{(\pi_1, \ldots, \pi_{\ell})})}.$$
Now, repeated use of the Murnaghan-Nakayama Rule (Theorem 3.4) shows that, for any $(\pi_1, \ldots, \pi_{\ell-1},\emptyset) \Vdash
w-|\pi_{\ell}|$, we have, writing $\pi_{\ell}=(k_1, \, \ldots , \, k_i)$, and since all irreducible characters of $L$ have degree 1,

$$\chi^{\alpha}({\cal D}_{(\pi_1, \ldots, \pi_{\ell})} )= \displaystyle
\sum_{\scriptstyle 1 \leq s_1, \ldots,s_i \leq \ell} \; \;
\sum_{\scriptstyle \tilde{\alpha} \in {\cal L}_{\alpha,(s_1,
\ldots,s_i)}^{\pi_{\ell}}} (-1)^{L_{\alpha \tilde{\alpha}}}
\chi^{\tilde{\alpha}}({\cal D}_{(\pi_1, \ldots,
\pi_{\ell-1})}^{w-|\pi_{\ell}|} ),$$where ${\cal L}_{\alpha,(s_1,
\ldots,s_i)}^{\pi_{\ell}}$ is the set of $\ell$-tuples of partitions
of $w-|\pi_{\ell}|$ which can be obtained from $\alpha$ by removing
successively a $k_1$-hook from $\alpha^{s_1}$, then a $k_2$-hook
from the ``$s_2$-th coordinate'' of the resulting $\ell$-tuple of
partitions of $w-k_1$, etc, and finally a $k_i$-hook from the
``$s_i$-th coordinate'' of the resulting $\ell$-tuple of partitions
of $w-(k_1+ \cdots +k_{i-1})$, and, for $ \tilde{\alpha} \in {\cal
L}_{\alpha,(s_1, \ldots,s_i)}^{\pi_{\ell}}$, $L_{\alpha
\tilde{\alpha}}$ is the sum of the leg-lengths of the hooks removed
to get from $\alpha$ to $\tilde{\alpha}$.

We therefore get, writing ${\bf s}=(s_1, \, \ldots , \, s_i)$,
$$\langle \chi^{\alpha}, \, \chi^{\beta}\rangle _{\scriptstyle
{\cal D}_{\pi_{\ell}}}= \displaystyle \frac{1}{d_{\pi_{\ell}}}
\sum_{\scriptstyle 1 \leq {\bf s} \leq \ell \atop 1 \leq {\bf t}
\leq \ell} \sum_{\scriptstyle \tilde{\alpha} \in {\cal
L}_{\alpha,{\bf s}}^{\pi_{\ell}} \atop \tilde{\beta} \in {\cal
L}_{\beta,{\bf t}}^{\pi_{\ell}} } (-1)^{L_{\alpha
\tilde{\alpha}}}(-1)^{L_{\beta \tilde{\beta}}} \langle
\chi^{\tilde{\alpha}}, \, \chi^{\tilde{\beta}} \rangle_{ {\cal
D}_{\emptyset}^{w- |\pi_{\ell}|}}.$$ We will use this formula to
exhibit a generalized perfect isometry between \\Irr$(L \wr S_w)$
and the principal $\ell$-block of $N \wr S_w$.

\subsection{The wreath product $N \wr S_w$}
We now turn to ${\cal N}=N_{S_{\ell w}}({\cal L})=N \wr S_w$. In particular, we want to find its principal $\ell$-block. Take $\{ g_1=\omega, \, g_2, \, \ldots, \, g_r \}$ representatives for the conjugacy classes of $N=N_{S_{\ell}}(L)$, where $\omega$ is an $\ell$-cycle generating $L$. For any $(\pi_1, \, \ldots , \, \pi_r) \Vdash w$, we write ${\cal S}_{(\pi_1, \ldots ,  \pi_r)} $ for the conjugacy class of elements of $N \wr S_w$ of cycle type $(\pi_1, \, \ldots , \, \pi_r)$.

\begin{defn} An element of $N \wr S_w$ of cycle type $(\pi_1, \, \ldots , \, \pi_r) \Vdash w$ is called $\ell${\emph{-regular}} if $\pi_1=\emptyset$, and $\ell${\emph{-singular}} otherwise.

We define the $\ell${\emph{-blocks}} of $N \wr S_w$ by orthogonality across the set ${\cal S}_{\emptyset}$ of $\ell$-regular elements.
\end{defn}

Note that, if $\ell$ is a prime, then we obtain the same $\ell$-regular elements and $\ell$-blocks if we identify $N \wr S_w$ with its natural imbedding in $S_{\ell w}$ and define an $\ell$-regular element according to its cyclic decomposition in $S_{\ell w}$. If the cyclic decomposition of $g_i$ in $S_{\ell}$ is $(m_1, \, , \ldots , \, m_s)$, then, for each $1 \leq j \leq s$, each $k$-cycle of $\pi_i$ gives $gcd(m_j,k)$ cycles of length $lcm(m_j,k)$. A cycle of length divisible by $\ell$ can only appear if $\ell$ divides $lcm(m_j,k)$. Since $k \leq w < \ell$, this forces $m_j$ to be divisible by $\ell$, and $g_i$ to be the $\ell$-cycle $g_1$ (in particular, the $\ell$-blocks are the ``ordinary'' prime blocks). This fails when $\ell$ is no longer prime.

\medskip
 We
partition $N \wr S_w$ into unions of conjugacy classes according to
$\pi_1 \in {\cal P}_{\leq w}$. We have
$$N \wr S_w = \displaystyle \coprod_{\pi_{1} \in {\cal P}_{\leq
w}} {\cal S}_{\pi_{1}},$$ where ${\cal S}_{\pi_{1}}$ is the
set of elements whose cycle type has $\pi_{1}$ as 1st part.

We write Irr$(N)=\{ \psi_1, \, \ldots , \, \psi_r \}$, and $\beta_0=\{ \psi_1, \, \ldots , \, \psi_{\ell}=1_{N} \}$ the principal $\ell$-block of $N$.

Given two $r$-tuples $\alpha$ and $\beta$ of partitions of $w$, we then have
$$\langle \chi^{\alpha}, \,
\chi^{\beta}\rangle _{\scriptstyle \ell-sing \atop N \wr S_w}=
\displaystyle \sum_{ \emptyset \neq \pi_1 \in {\cal P}_{\leq
w}} \langle \chi^{\alpha}, \, \chi^{\beta}\rangle _{\scriptstyle
{\cal S}_{\pi_1}}.$$ We have, for any $\emptyset \neq
\pi_1 \in {\cal P}_{\leq w}$,
$$\langle \chi^{\alpha}, \, \chi^{\beta}\rangle _{\scriptstyle
{\cal S}_{\pi_1}}= \displaystyle \frac{1}{|N \wr S_w|}
\sum_{\scriptstyle (\emptyset, \pi_2, \ldots, \pi_{r}) \Vdash
w-|\pi_{1}|} |{\cal S}_{(\pi_1, \ldots, \pi_{r})}|
\chi^{\alpha}({\cal S}_{(\pi_1, \ldots, \pi_{r})} )
\overline{\chi^{\beta}( {\cal S}_{(\pi_1, \ldots, \pi_r)})}.$$

Now, as in the previous section, for any $(\emptyset, \pi_2, \ldots, \pi_r) \Vdash
w-|\pi_1|$, Lemma 3.2 gives, writing $b_{jk}$ for the number of $k$-cycles in $\pi_j$
($1\leq j \leq \ell$, $1 \leq k \leq w$),
$$\begin{array}{rl} \displaystyle \frac{|N \wr
S_w|}{|{\cal S}_{(\pi_1, \ldots, \pi_r)}|} & =   \displaystyle \prod_{\scriptstyle 1 \leq k \leq w} b_{1 k}! (k |C_N(g_1)|)^{b_{1 k}} \displaystyle \prod_{\scriptstyle 2\leq j \leq r-1 \atop 1 \leq k \leq w} b_{jk}! (k |C_N(g_j)|)^{b_{jk}} \\ & = \displaystyle \left( \prod_{\scriptstyle 1 \leq k \leq w} b_{1 k}! (k \ell)^{b_{\ell k}} \right) \displaystyle \frac{|N \wr
S_{w-|\pi_1|}|}{|{\cal S}_{(\emptyset, \pi_2, \ldots, \pi_r)}^{w- |\pi_1|}|}, \end{array}
$$
where ${\cal S}_{(\emptyset , \pi_2, \ldots, \pi_r)}^{w- |\pi_1|}$ denotes the set of elements of cycle type $(\emptyset, \pi_2, \ldots, \pi_r)$ in $N \wr
S_{w-|\pi_1|}$. Writing $c_{\pi_1}=\prod_{\scriptstyle 1 \leq k \leq w} b_{1 k}! (k \ell)^{b_{1 k}}$, we get

$$\langle \chi^{\alpha}, \, \chi^{\beta}\rangle _{\scriptstyle
{\cal S}_{\pi_1}}= \displaystyle \frac{1}{c_{\pi_1}}
\sum_{\scriptstyle (\emptyset,\pi_2, \ldots, \pi_r) \Vdash
w-|\pi_1|}  \frac{|{\cal S}_{(\emptyset, \pi_2, \ldots, \pi_r)}^{w- |\pi_r|}|}{|N \wr
S_{w-|\pi_1|}|}
\chi^{\alpha}({\cal S}_{(\pi_1, \ldots, \pi_r)} )
\overline{\chi^{\beta}( {\cal S}_{(\pi_1, \ldots, \pi_r)})}.$$
Now, as before, if we write $\pi_{1}=(k_1, \, \ldots , \, k_i)$, then repeated use the Murnaghan-Nakayama Rule (Theorem 3.4) shows that, for any $(\emptyset, \pi_2, \ldots, \pi_r) \Vdash
w-|\pi_1|$, we have, with the notations we introduced,

$$\chi^{\alpha}({\cal S}_{(\pi_1, \ldots, \pi_r)} )= \displaystyle
\sum_{\scriptstyle 1 \leq {\bf s} \leq r} \sum_{\scriptstyle
\tilde{\alpha} \in {\cal L}_{\alpha,{\bf s}}^{\pi_1}}
 \Psi_{{\bf s}}(g_1) (-1)^{L_{\alpha \tilde{\alpha}}}
 \chi^{\tilde{\alpha}}({\cal S}_{(\pi_2, \ldots, \pi_{r})}^{w-|\pi_{1}|} ).$$
where $\Psi_{{\bf s}}(g_1)=\psi_{s_1}(g_1) \ldots \psi_{s_i}(g_1)$.
We therefore get
$$\langle \chi^{\alpha}, \, \chi^{\beta}\rangle _{\scriptstyle
{\cal S}_{\pi_1}}= \displaystyle \frac{1}{c_{\pi_1}}
\sum_{\scriptstyle 1 \leq {\bf s} \leq r \atop 1 \leq {\bf t} \leq
r}  \; \sum_{\scriptstyle \tilde{\alpha} \in {\cal L}_{\alpha,{\bf
s}}^{\pi_1}
 \atop \tilde{\beta} \in {\cal L}_{\beta,{\bf s}}^{\pi_1} }  \Psi_{{\bf s}}(g_1)
   \overline{\Psi_{{\bf t}}(g_1)} (-1)^{L_{\alpha \tilde{\alpha}}}(-1)^{L_{\beta \tilde{\beta}}}
    \langle \chi^{\tilde{\alpha}}, \, \chi^{\tilde{\beta}} \rangle_{ {\cal S}_{\emptyset}^{w- |\pi_1|}}.
    \; \; (\dagger)$$
Note that, if $\pi_1=\pi_{\ell}$, then $d_{\pi_{\ell}}=c_{\pi_1}$. This will allow us later to exhibit a generalized perfect isometry.

We first prove the following:

\begin{thm}
If, for some $\alpha , \beta \Vdash w$, $\alpha \neq \beta$, and
$\pi_1 \in {\cal P}_{\leq w}$, we have $\langle \chi^{\alpha}, \,
\chi^{\beta}\rangle _{\scriptstyle {\cal S}_{\pi_1}} \neq 0$, then,
for any $\ell < k \leq r$, we have $\alpha^k=\beta^k=\emptyset$. In
particular, the principal $\ell$-block of $N \wr S_w$ is included in
$\{ \chi^{\alpha}, \, \alpha=(\alpha_1, \ldots, \alpha_{\ell},
\emptyset, \ldots, \emptyset) \Vdash w \}$ and, if $\alpha^k \neq
\emptyset$ for some $\ell <k \leq r$, then $\{ \chi^{\alpha} \}$ is
an $\ell$-block of $N \wr S_w$.
\end{thm}

\begin{proof}
We use induction on $w$. For $w=1$, the result is just the description of the $\ell$-blocks of $N$. Indeed, any element of $N \wr S_w$ has cycle type $(\pi, \ldots, \pi_r)$ where $\pi_i=(1)$ for some $1 \leq i \leq r$ and $\pi_j= \emptyset$ for $j \neq i$ (in particular, it is $\ell$-regular if and only if $i \neq 1$ (i.e. $\pi_1 = \emptyset$) and $\ell$-singular if and only if $i=1$ (i.e. $\pi_1=(1)$)). For any $\alpha \Vdash w$, we have $\alpha^i=(1)$ for some $1 \leq i \leq r$ and $\alpha^j= \emptyset$ for $j \neq i$, and $\chi^{\alpha}=\psi_i$.

\medskip

Now take $w >1$, and suppose the result true up to $w-1$. Suppose $\langle \chi^{\alpha}, \, \chi^{\beta}\rangle _{\scriptstyle
{\cal S}_{\pi_1}} \neq 0$ for some $\pi_1 \in {\cal P}_{\leq w}$.

\smallskip

Suppose first that $\pi_1 \neq \emptyset$. Then, by $(\dagger)$, if
we write $\pi_1=(k_1, \ldots, k_i)$, then there exist $1 \leq s_1,
\ldots , s_i, t_1, \ldots , t_i \leq r$, there exist $\tilde{\alpha}
\in {\cal L}_{\alpha,{\bf s}}^{\pi_1}$ and $ \tilde{\beta} \in {\cal
L}_{\beta,{\bf t}}^{\pi_1}$ such that
$$\psi_{s_1}(g_1)  \overline{\psi_{t_1}(g_1)} \ldots \psi_{s_i}(g_1)
\overline{\psi_{t_i}(g_1)} \langle \chi^{\tilde{\alpha}}, \,
\chi^{\tilde{\beta}} \rangle_{ {\cal S}_{\emptyset}^{w- |\pi_1|}}
\neq 0.$$ By the study of the $\ell$-blocks of $N$, we know that
$\psi_k(g_1) \neq 0$ if and only if $1 \leq k \leq \ell$. Thus $1
\leq s_1, \ldots , s_i, t_1, \ldots , t_i \leq \ell$, so that, if
$k> \ell$, then $\alpha^k=\tilde{\alpha}^k$ and $\beta^k =
\tilde{\beta}^k$ (since only hooks from $\alpha^{s_1}, \, \ldots, \,
\alpha^{s_i}$ are removed to get from $\alpha$ to $\tilde{\alpha}$,
and similarly for $\beta$ and $\tilde{\beta}$). And, by induction
hypothesis (since $w- |\pi_1| < w$), $\langle \chi^{\tilde{\alpha}},
\, \chi^{\tilde{\beta}} \rangle_{ {\cal S}_{\emptyset}^{w- |\pi_1|}}
\neq 0$ implies that $\tilde{\alpha}^k=\tilde{\beta}^k=\emptyset$
for $k > \ell$. Hence $\alpha^k=\beta^k=\emptyset$ for $k > \ell$.

\smallskip
Next suppose that $\pi_1 = \emptyset$. Then, since $N \wr S_w = \coprod_{\pi_1 \in {\cal P}_{\leq w}} {\cal S}_{\pi_1}$, we have
$$0=\langle \chi^{\alpha}, \, \chi^{\beta} \rangle  = \displaystyle \sum_{\pi_1 \in {\cal P}_{\leq w}} \langle \chi^{\alpha}, \, \chi^{\beta}\rangle _{\scriptstyle
{\cal S}_{\pi_1}},$$
so that there exists $\emptyset \neq \pi \in {\cal P}_{\leq w}$ such that $\langle \chi^{\alpha}, \, \chi^{\beta}\rangle _{\scriptstyle
{\cal S}_{\pi}} \neq 0$. By the previous case, this implies that $\alpha^k=\beta^k=\emptyset$ for $k > \ell$. This proves our assertion.

\end{proof}

We now come to the construction of a generalized perfect isometry
between the principal $\ell$-block of $N \wr S_w$ and Irr$(L \wr
S_w)$. In order to relate the contribution of two characters of $N
\wr S_w$ on $\ell$-singular elements to an inner product on singular
elements in $L \wr S_w$, we want the $\psi_i's$ which take value
$-1$ on the $\ell$-cycle $g_1$ to ``disappear'' in the
Murnaghan-Nakayama Rule. We achieve this by using the following
transformation. Let $\{ \psi_1, \, \ldots , \, \psi_m \}$ be all the
irreducible characters of $N$ which take value $-1$ on $g_1$. For
any $\alpha = (\alpha^1, \ldots, \alpha^{\ell},  \emptyset, \ldots,
\emptyset) \Vdash w$, we let $\alpha^*=(\overline{\alpha^1}, \ldots,
\overline{\alpha^m}, \alpha^{m+1} , \ldots , \alpha^{\ell},
\emptyset, \ldots, \emptyset) \Vdash w$, where $^{\overline{ \mbox{
} }}$ denotes conjugation of partitions (i.e. multiplication of the
corresponding character of the symmetric group by the signature) and
we set
$$\chi^{\alpha}_0 = \displaystyle (-1)^{\sum_{n=1}^m | \alpha^n |} \chi^{\alpha^*}.$$
Note that $^*$ is a bijection from $\{ \alpha=(\alpha_1, \ldots, \alpha_{\ell}, \emptyset, \ldots, \emptyset) \Vdash w \}$ onto itself.

\begin{prop}
Take $\alpha = (\alpha^1, \ldots, \alpha^{\ell},  \emptyset, \ldots, \emptyset) \Vdash w$, and take $g
\in N \wr S_w$ of cycle type $(a_{ij}(g))_{1 \leq i \leq r, 1 \leq j
\leq w}$ such that $a_{1k}(g) >0$ for some $1
\leq k \leq w$. Let $\rho \in N \wr S_{w-k}$ be of cycle type
$(a_{ij}(\rho))_{1 \leq i \leq r, 1 \leq j \leq w}$, where
$$(a_{ij}(\rho)) = \left\{ \begin{array}{ll} a_{1k}(g)-1 &
{\mbox{if}} \; i=1 \; \mbox{and} \; j=k \\ a_{ij}(g) &
\mbox{otherwise} \end{array} \right. $$
Then we have, with the notations of Theorem 3.4,
$$\chi^{\alpha}_0(g)= \displaystyle \sum_{s=1}^{\ell}
\sum_{h_{ij}^{\alpha^s}=k} (-1)^{L_{ij}^{\alpha^s}}
\chi^{\alpha-R_{ij}^{\alpha^s}}_0(\rho).$$
\end{prop}

\begin{proof}
By Theorem 3.4, we have
$$\begin{array}{rl} \chi^{\alpha}_0(g) & = \displaystyle (-1)^{\sum_{n=1}^m | \alpha^n |} \chi^{\alpha^*}(g) \\ & =\displaystyle (-1)^{\sum_{n=1}^m | \alpha^n |}   \sum_{s=1}^{r} \psi_s(g_1)
\sum_{h_{ij}^{\alpha^{*s}}=k} (-1)^{L_{ij}^{\alpha^{*s}}}
\chi^{\alpha^*-R_{ij}^{\alpha^{*s}}}(\rho) \\ & = \displaystyle  \sum_{s=1}^{\ell} \chi^{\alpha}_0(g)_s \end{array}$$
(since, as in the proof of 3.7, if $s \geq \ell$, then $\chi^{\alpha}_0(g)_s=\psi_s(g_1)=0$).

\smallskip
Suppose first that $m+1 \leq s \leq \ell$. Then $\psi_s(g_1)=1$ and $\alpha^{*s}=\alpha^s$. Thus
$$\psi_s(g_1)
\sum_{h_{ij}^{\alpha^{*s}}=k} (-1)^{L_{ij}^{\alpha^{*s}}}
\chi^{\alpha^*-R_{ij}^{\alpha^{*s}}}(\rho)=
\sum_{h_{ij}^{\alpha^s}=k} (-1)^{L_{ij}^{\alpha^s}}
\chi^{\alpha^*-R_{ij}^{\alpha^s}}(\rho),$$ and
$$\chi^{\alpha^*-R_{ij}^{\alpha^s}}=\chi^{(\alpha-R_{ij}^{\alpha^s})^*}=
\displaystyle (-1)^{\sum_{n=1}^m | (\alpha-R_{ij}^{\alpha^s})^n |}
\chi^{\alpha-R_{ij}^{\alpha^s}}_0 = \displaystyle (-1)^{\sum_{n=1}^m
| \alpha^n |} \chi^{\alpha-R_{ij}^{\alpha^s}}_0,$$ so that
$$\chi^{\alpha}_0(g)_s = \displaystyle \sum_{h_{ij}^{\alpha^s}=k}
 (-1)^{L_{ij}^{\alpha^s}} \chi^{(\alpha-R_{ij}^{\alpha^s})}_0(\rho).$$

Now suppose $1 \leq s \leq m$. Then $\psi_s(g_1)=-1$. The set of hook-lengths of any given partition is preserved by conjugation, and, for any hook-length $k$, conjugation gives a bijection from the set of $k$-hooks of $\alpha^{*s}$ onto the set of $k$-hooks of $\alpha^s$. For any $h_{ij}^{\alpha^{*s}}=k$, we have $L_{ij}^{\alpha^{*s}}=k-L_{ij}^{\alpha^s}-1$ and $\chi^{\alpha^*-R_{ij}^{\alpha^{*s}}}=\chi^{(\alpha-R_{ij}^{\alpha^s})^*}$. Thus
$$\begin{array}{rl} \chi^{\alpha}_0(g)_s & = \displaystyle (-1)^{\sum_{n=1}^m | \alpha^n |}  \psi_s(g_1)
\sum_{h_{ij}^{\alpha^{*s}}=k} (-1)^{L_{ij}^{\alpha^{*s}}}
\chi^{\alpha^*-R_{ij}^{\alpha^{*s}}}(\rho) \\ & = \displaystyle
\sum_{h_{ij}^{\alpha^s}=k} (-1)^{\sum_{n=1}^m | \alpha^n |} (-1)
(-1)^{k-L_{ij}^{\alpha^s}-1}
\chi^{(\alpha-R_{ij}^{\alpha^s})^*}(\rho) \\ & = \displaystyle
\sum_{h_{ij}^{\alpha^s}=k} (-1)^{\sum_{n=1}^m |
(\alpha-R_{ij}^{\alpha^s})^n |+k} (-1)^{k-L_{ij}^{\alpha^s}}
\chi^{(\alpha-R_{ij}^{\alpha^s})^*}(\rho) \\ &= \displaystyle
\sum_{h_{ij}^{\alpha^s}=k}
 (-1)^{L_{ij}^{\alpha^s}} \chi^{(\alpha-R_{ij}^{\alpha^s})}_0(\rho). \end{array}$$
This concludes the proof.

\end{proof}

\noindent
{\bf{Remark:}} Of course, the result of this proposition is precisely the motivation for the introduction of $^*$. We want the virtual character $\chi^{\alpha}_0$ to satisfy a convenient analogue of the Murnaghan-Nakayama Rule when evaluated at $\ell$-singular elements. If we go through the proof of Theorem 4.4 in \cite{Pfeiffer}, we see that, essentially, our equality holds because of the following general fact (which one proves easily)

\begin{prop}
Let $n \geq 1$ be any integer, and, for $\pi \in S_n$, let $c(\pi)$ be the number of cycles (including the trivial ones) in $\pi$. Define
$$\eta \colon \left\{ \begin{array}{rl} S_n \longrightarrow & {\bf C} \\ \pi \longmapsto & (-1)^{c(\pi)} \end{array} \right.$$
Then $\eta = (-1)^n \varepsilon_n$, where $\varepsilon_n$ is the signature of $S_n$.
\end{prop}

Now, take any $\pi_{1}=(k_1, \, \ldots , \, k_i) \in {\cal P}_{\leq
w}$. Then repeated use of Proposition 3.8 yields, for any
$\alpha, \beta \Vdash w$, the following analogue of $(\dagger)$

$$\langle \chi^{\alpha}_0, \, \chi^{\beta}_0\rangle _{\scriptstyle
{\cal S}_{\pi_1}}= \displaystyle \frac{1}{c_{\pi_1}}
\sum_{\scriptstyle 1 \leq {\bf s} \leq \ell \atop 1 \leq {\bf t}
\leq \ell} \sum_{\scriptstyle \tilde{\alpha} \in {\cal
L}_{\alpha,{\bf s}}^{\pi_1} \atop \tilde{\beta} \in {\cal
L}_{\beta,{\bf t}}^{\pi_1} } (-1)^{L_{\alpha
\tilde{\alpha}}}(-1)^{L_{\beta \tilde{\beta}}} \langle
\chi^{\tilde{\alpha}}_0, \, \chi^{\tilde{\beta}}_0 \rangle_{ {\cal
S}_{\emptyset}^{w- |\pi_1|}}. \; \; \; \; \; \; (\dagger \dagger)$$

This equality provides the main ingredient in our generalized
perfect isometry:

\begin{prop}
Take any $\alpha=(\alpha_1, \, \ldots , \, \alpha_{\ell}) \Vdash w$
and $\beta=(\beta_1, \, \ldots , \, \alpha_{\ell}) \Vdash w$. Define
two $r$-tuples of partitions $\alpha_r=(\alpha_1, \, \ldots , \,
\alpha_{\ell}, \, \emptyset, \, \ldots , \, \emptyset) $ and\\
$\beta_r=(\beta_1, \, \ldots , \, \alpha_{\ell}, \, \emptyset , \,
\ldots , \, \emptyset)$ of $ w$. Then
$$\langle \chi^{\alpha_r}_0, \, \chi^{\beta_r}_0\rangle _{\scriptstyle
\ell-reg \atop N \wr S_w} = \langle \chi^{\alpha}, \,
\chi^{\beta}\rangle _{\scriptstyle reg \atop L \wr S_w}.$$

\end{prop}

\begin{proof}

By the results of section 3.2, we have
$$ \begin{array}{rl} \langle \chi^{\alpha}, \,
\chi^{\beta}\rangle _{\scriptstyle sing \atop L \wr S_w} & =
\displaystyle \sum_{ \emptyset \neq \pi_{\ell} \in {\cal P}_{\leq
w}} \frac{1}{d_{\pi_{\ell}}} \sum_{\scriptstyle 1 \leq {\bf s}, {\bf
t}  \leq \ell} \sum_{\scriptstyle \tilde{\alpha} \in {\cal
L}_{\alpha, {\bf s}}^{\pi_{\ell}} \atop \tilde{\beta} \in {\cal
L}_{\beta, {\bf t}}^{\pi_{\ell}} } (-1)^{L_{\alpha
\tilde{\alpha}}}(-1)^{L_{\beta \tilde{\beta}}} \langle
\chi^{\tilde{\alpha}}, \, \chi^{\tilde{\beta}} \rangle_{ {\cal
D}_{\emptyset}^{w- |\pi_{\ell}|}} \\ & = \displaystyle \sum_{k=1}^w
 \sum_{  \pi_{\ell} \vdash k} \frac{1}{d_{\pi_{\ell}}} \sum_{\scriptstyle 1 \leq {\bf s}, {\bf
t}  \leq \ell} \sum_{\scriptstyle \tilde{\alpha} \in {\cal
L}_{\alpha, {\bf s}}^{\pi_{\ell}} \atop \tilde{\beta} \in {\cal
L}_{\beta, {\bf t}}^{\pi_{\ell}} } (-1)^{L_{\alpha
\tilde{\alpha}}}(-1)^{L_{\beta \tilde{\beta}}} \langle
\chi^{\tilde{\alpha}}, \, \chi^{\tilde{\beta}} \rangle_{ {\cal
D}_{\emptyset}^{w- k}}.
\end{array}$$
On the other hand, by $(\dagger \dagger)$, we have
$$\langle \chi^{\alpha_r}_0, \, \chi^{\beta_r}_0\rangle _{\scriptstyle
\ell-reg \atop N \wr S_w} = \displaystyle \sum_{k=1}^w
 \sum_{  \pi_{1} \vdash k} \frac{1}{c_{\pi_{1}}} \sum_{1 \leq {\bf s},  {\bf t} \leq \ell }
 \sum_{\scriptstyle \tilde{\alpha} \in {\cal
L}_{\alpha, {\bf s}}^{\pi_1} \atop \tilde{\beta} \in {\cal
L}_{\beta, {\bf t}}^{\pi_1} }  (-1)^{L_{\alpha_r
\tilde{\alpha}_r}}(-1)^{L_{\beta_r \tilde{\beta}_r}} \langle
\chi^{\tilde{\alpha}_r}, \, \chi^{\tilde{\beta}_r} \rangle_{ {\cal
S}_{\emptyset}^{w- k}}.$$ We can now easily prove the result by
induction on $w$. By the previous two equalities, if the result is
true up to $w-1$, then it is true for $w$ ($w \geq 1$). Finally, the
result is clearly true in the case $w=0$, since ${\cal
S}_{\emptyset}^0={\cal D}_{\emptyset}^0=\emptyset$, and
$\chi^{\emptyset}=\chi^{\emptyset_r}=1$.
\end{proof}

In particular, we have the following

\begin{cor}
The principal $\ell$-block of $N \wr S_w$ is $$b_0=\{
\chi^{\alpha}_r, \, \alpha_r=(\alpha_1, \ldots, \alpha_{\ell},
\emptyset, \ldots, \emptyset) \Vdash w \}.$$

\end{cor}

\section{A generalized perfect isometry}

We can now prove our main result:

\begin{thm}

Let $\ell \geq 2$ be any integer. Let $G$ be the symmetric group
$S_{\ell w + r}$, where $0 \leq r, \, w < \ell$, and let ${\cal L}
\cong L^w \cong {\bf Z}_{\ell}^w$ be a "Sylow $\ell$-subgroup" of
$G$. Then there is a generalized perfect isometry with respect to
$\ell$-regular elements between the principal $\ell$-blocks of $G$
and $N_G({\cal L})$.
\end{thm}

\begin{proof}
The principal $\ell$-block of $S_{\ell w + r}$ is $B_0= \{
\varphi_{\lambda}, \, \lambda \in {\cal P}_r \}$, where ${\cal P}_r$
denotes the set of partitions of $\ell w +r$ with $\ell$-core $(r)$.
For any $\lambda \in {\cal P}_r$, we write $$\alpha_{\lambda}=(
\alpha_{\lambda}^1, \, \ldots , \, \alpha_{\lambda}^{\ell}) \Vdash
w$$ the $\ell$-quotient of $\lambda$, and, as before,
$$\alpha_{\lambda , r}=( \alpha_{\lambda}^1, \, \ldots , \,
\alpha_{\lambda}^{\ell}, \, \emptyset , \, \ldots , \, \emptyset)
\Vdash w$$and$$\alpha_{\lambda , r}^*=(
\overline{\alpha_{\lambda}^1}, \, \ldots , \,
\overline{\alpha_{\lambda}^m}, \, \alpha_{\lambda}^{m+1}, \, \ldots
, \, \alpha_{\lambda}^{\ell}, \, \emptyset , \, \ldots , \,
\emptyset) \Vdash w.$$ We have Irr$(L \wr S_w) = \{
\chi^{\alpha_{\lambda}}, \, \lambda \in {\cal P}_r \}$, and, by
Theorem 3.5, there exist signs $\{ \varepsilon(\lambda) , \, \lambda
\in {\cal P}_r \}$ such that, for any $\lambda, \mu \in {\cal P}_r$,
$$\langle \varphi_{\lambda}, \, \varphi_{\mu} \rangle _{\scriptstyle
\ell-reg \atop S_{\ell w +r}} =\langle \varepsilon (\lambda)
\chi^{\alpha_{\lambda}}, \, \varepsilon (\mu)
\chi^{\alpha_{\mu}}\rangle _{\scriptstyle reg \atop L \wr S_w}.$$
For any $\lambda \in {\cal P}_r$, we write $\eta( \lambda) =
(-1)^{\sum_{n=1}^m | \alpha_{\lambda}^n |}$ (so that
$\chi^{\alpha_{\lambda , r}}_0= \eta( \lambda) \chi^{\alpha_{\lambda
, r}^*}$). Then $\chi^{\alpha_{\lambda}} \longmapsto
\chi^{\alpha_{\lambda , r}^*}$ is a bijection from Irr$(L \wr S_w)$
onto the principal $\ell$-block $b_0= \{ \chi^{\alpha_{\lambda ,
r}^*}, \, \lambda \in {\cal P}_r \}$ of $N \wr S_w$, and, for any
$\lambda, \, \mu \in {\cal P}_r$, we have, by Proposition 3.10,

$$\langle \chi^{\alpha_{\lambda}}, \, \chi^{\alpha_{\mu}}\rangle
_{\scriptstyle reg \atop L \wr S_w}= \langle \eta(\lambda)
\chi^{\alpha_{\lambda , r}^*}, \, \eta(\mu) \chi^{\alpha_{\mu ,
r}^*}\rangle _{\scriptstyle \ell-reg \atop N \wr S_w}.$$ Finally,
using the results of section 2.1, we get that
$$I \colon \left\{ \begin{array}{ccl} B_0 & \longrightarrow & b_0
\otimes \{1_{S_r} \} \\ \varphi_{\lambda} & \longmapsto &
\chi^{\alpha_{\lambda , r}^*} \otimes 1_{S_r} \end{array} \right.$$
is a bijection between the principal $\ell$-blocks of $G=S_{\ell
w+r}$ and $N_G({\cal L})$ such that, for any $\lambda, \, \mu \in
{\cal P}_r$, we have
$$\langle \varphi_{\lambda}, \, \varphi_{\mu} \rangle _{\scriptstyle
\ell-reg \atop S_{\ell w +r}} = \langle \varepsilon (\lambda) \eta
(\lambda) I( \varphi_{\lambda}), \, \varepsilon (\mu) \eta (\mu) I
(\varphi_{\mu}) \rangle _{\scriptstyle \ell-reg \atop N_G({\cal
L})},$$which ends the proof.
\end{proof}

\bigskip
\noindent Many thanks to Elvis.

\end{document}